\documentclass{amsart}

\usepackage[latin1]{inputenc}
\usepackage{amsfonts}
\usepackage{amssymb}
\usepackage{eufrak}
\usepackage{amsmath}
\usepackage{amsthm}
\usepackage{a4wide}
\usepackage[all]{xy}
\usepackage{textcomp}

\newtheorem{thm}{Theorem}[section]
\newtheorem{defi}[thm]{Definition}

\newtheorem{lemma}[thm]{Lemma}

\begin{document}

\title[On iterated translated points]
{On iterated translated points for contactomorphisms of $\mathbb{R}^{2n+1}$ and $\mathbb{R}^{2n}\times S^1$}

\author{Sheila Sandon}
\address{Laboratoire de Math\'{e}matiques Jean Leray, Universit\'{e} de Nantes, 2, rue de la Houssini\`{e}re, 44322 Nantes, France}
\email{sheila.sandon@univ-nantes.fr}

\begin{abstract}
\noindent
A point $q$ in a contact manifold is called a \textit{translated point} for a contactomorphism $\phi$ with respect to some fixed contact form if $\phi(q)$ and $q$ belong to the same Reeb orbit and the contact form is preserved at $q$. The problem of existence of translated points has an interpretation in terms of Reeb chords between Legendrian submanifolds, and can be seen as a special case of the problem of leafwise coisotropic intersections. For a compactly supported contactomorphism $\phi$ of $\mathbb{R}^{2n+1}$ or $\mathbb{R}^{2n}\times S^1$ contact isotopic to the identity, existence of translated points follows immediately from Chekanov's theorem on critical points of quasi-functions and Bhupal's graph construction. In this article we prove that if $\phi$ is positive then there are infinitely many non-trivial geometrically distinct iterated translated points, i.e. translated points of some iteration $\phi^k$. This result can be seen as  a (partial) contact analogue of the result of Viterbo on existence of infinitely many iterated fixed points for compactly supported Hamiltonian symplectomorphisms of $\mathbb{R}^{2n}$, and is obtained with generating functions techniques in the setting of \cite{mio1}.
\end{abstract}

\maketitle

\section{Introduction}\label{intro}

Let $\phi$ be a Hamiltonian symplectomorphism of a closed symplectic manifold $(W,\omega)$. Arnold's famous conjecture states that $\phi$ should have at least as many fixed points as a function on $W$ must have critical points. This conjecture was posed in the 60's and is still not proved in full generality. It has been one of the main motivating forces in the development of some of the most powerful techniques in modern symplectic topology. Consider now a contact manifold $\big(M,\xi=\text{ker}(\alpha)\big)$ and let $\phi$ be a contactomorphism, contact isotopic to the identity. In contrast to the case of a Hamiltonian symplectomorphism, $\phi$ does not need to have any fixed point. We can consider for example the Reeb flow on $M$, i.e. the contact isotopy generated by the vector field $R_{\alpha}$ defined by the conditions $\iota_{R_{\alpha}}d\alpha=0$ and $\alpha(R_{\alpha})=1$. Since the Reeb vector field never vanishes, for small $t$ the time-$t$ map of the Reeb flow does not have any fixed point. This difference with respect to the symplectic case is reflected by the fact that while in the symplectic case the Hamiltonian isotopy generated by a constant Hamiltonian function is constantly the identity, in the contact case a constant Hamiltonian generates the Reeb flow. Motivated by this, we can consider instead an analogue of the Arnold conjecture for \textit{translated points} of $\phi$.

\begin{defi}
Let $\phi$ be a contactomorphism of a (cooriented) contact manifold $\big(M,\xi=\text{ker}(\alpha)\big)$ and let $g:M\rightarrow \mathbb{R}$ be the function defined by $\phi^{\ast}\alpha=e^g\alpha$. A point $q$ of $M$ is called a \textbf{translated point} for $\phi$ (with respect to the contact form $\alpha$) if $\phi(q)$ and $q$ belong to the same Reeb orbit and $g(q)=0$.
\end{defi}

Consider for example the prequantization of an exact symplectic manifold $(M,\omega=-d\lambda)$, i.e. the manifold $M\times\mathbb{R}$ endowed with the contact structure $\text{ker}(d\theta-\lambda)$ where $\theta$ is the $\mathbb{R}$-coordinate. Then every point in the fiber above a fixed point of a Hamiltonian symplectomorphism $\phi$ of $M$ is a translated point of the lift $\widetilde{\phi}$. Recall that a lift $\widetilde{\phi}$ of $\phi$ is a contactomorphism of $M\times\mathbb{R}$ defined by $\widetilde{\phi}(x,\theta)=\big(\phi(x),\theta+F(x)\big)$ where $F$ is a function on $M$ satisfying $\phi^{\ast}\lambda-\lambda=dF$. Note that the condition on $g$ is automatically satisfied in this case because $\widetilde{\phi}$ preserves the contact form everywhere.\\
\\
Our motivation in considering the notion of translated points comes not so much from the previous example but rather from the special case of compactly supported contactomorphisms of $M=\big(\mathbb{R}^{2n}\times S^1,\xi=\text{ker}\,(dz-ydx)\big)$. Here translated points correspond to critical points of generating functions, and are at the base of the constructions of a contact capacity and contact homology groups for domains and of a bi-invariant partial order and a bi-invariant metric on the contactomorphism group \cite{B,mio1,mio2}.

\begin{defi}
A point $q$ of $M$ is called an \textbf{iterated translated point} for a contactomorphism $\phi$ if it is a translated point for some iteration $\phi^k$.
\end{defi}

We can consider the following problems.\\
\\
\textbf{PROBLEM 1}: What can we say about existence of translated points? Can we prove an analogue of the Arnold conjecture, i.e. that on a closed $M$ every contactomorphism which is contact isotopic to the identity must have at least as many translated points as a function on $M$ must have critical points? \\
\\
\textbf{PROBLEM 2}: Can we prove that on a closed manifold every contactomorphism which is contact isotopic to the identity has infinitely many iterated translated points? \\
\\
As we will now discuss, Problem 1 has a interpretation in terms of Reeb chords. Recall first that the Arnold conjecture for fixed points of Hamiltonian symplectomorphisms can be seen as a statement on Lagrangian intersections. Indeed, if $\phi$ is a Hamiltonian symplectomorphism of $(W,\omega)$ then its graph is a Lagrangian submanifold of the product $(W\times W,-\omega\oplus\omega)$, and the fixed points of $\phi$ correspond to intersections of the graph with the diagonal. A similar interpretation also exists for the problem of translated points, but Lagrangian intersections are replaced by Reeb chords between Legendrian submanifolds. Let $\phi$ be a contactomorphism of a contact manifold $\big(M,\xi=\text{ker}(\alpha)\big)$ and consider $M\times M\times \mathbb{R}$ with the product contact form $A=e^{\theta}\alpha_1-\alpha_2$. Here $\alpha_1$ and $\alpha_2$ are the pullback of $\alpha$ with respect to the projections of $M\times M\times \mathbb{R}$ into the first and second factors respectively, and $\theta$ is the $\mathbb{R}$-coordinate. Note that the Reeb flow of $A$ is generated by the vector field $(0,-R_{\alpha},0)$ where $R_{\alpha}$ is the Reeb vector field of $\alpha$. We define the diagonal of $M\times M\times \mathbb{R}$ to be the Legendrian submanifold $\{(q,q,0)\;|\; q\in M\}$, and the graph of a contactomorphism $\phi$ of $M$ to be the image of the Legendrian embedding $M \rightarrow M\times M\times \mathbb{R}$, $q \mapsto \big(q,\phi(q),g(q)\big)$ where $g$ is the function defined by $\phi^{\ast}\alpha=e^g\alpha$. If $q\in M$ is a translated point of $\phi$ then the corresponding point in the graph is of the form $\big(q,(\varphi^{R_{\alpha}})_{t_0}(q),0\big)$ for some $t_0$, where $(\varphi^{R_{\alpha}})_t$ is the Reeb flow of $\alpha$. Note that this point is in the same Reeb orbit as $(q,q,0)$, which belongs to the diagonal. Hence, translated points of $\phi$ correspond to Reeb chords between its graph and the diagonal in $M\times M\times \mathbb{R}$. The problem of counting translated points of the contactomorphism $\phi$ of $M$ is thus reduced to the problem of counting Reeb chords between its graph and the diagonal. Note that these are two Legendrian submanifolds that are contact isotopic to each other. The existence problem for Reeb chords connecting two Legendrian submanifolds contact isotopic to each other can be considered as a relative version of another famous conjecture by Arnold, the chord conjecture. This was stated in 1986 \cite{A} originally as the assertion that every Legendrian knot in the standard contact $S^3$ should have a Reeb chord. The generalized chord conjecture states existence of Reeb chords for Legendrian submanifolds of any closed contact manifold $(M^{2n-1},\xi)$. Partial results about the chord conjecture were obtained by Givental \cite{Giv}, Abbas \cite{Ab1,Ab2}, Mohnke \cite{Mo}, Cieliebak \cite{Ci} and Hutchings and Taubes \cite{HT}. As far as I know, the only result in the literature about existence of Reeb chords connecting two different Legendrian submanifolds is Chekanov's theorem \cite{C}, which concerns Reeb chords between the 0-section in a 1-jet bundle and a contact deformation of it. As we will now discuss, Chekanov's result implies existence of translated points in the special case of compactly supported contactomorphisms of $\mathbb{R}^{2n+1}$ and $\mathbb{R}^{2n}\times S^1$.\\
\\
Let $B$ be a smooth closed manifold, and consider in the 1-jet bundle $J^1B=T^{\ast}B\times\mathbb{R}$ a Legendrian submanifold $L$ which is contact isotopic to the 0-section. Note that the Reeb flow of $J^1B$ with respect to the standard contact form $d\theta-pdq$ is given by translation in the $\mathbb{R}$-direction. Thus Reeb chords connecting $L$ to the 0-section correspond to intersections of $L$ with the 0-wall, which is defined to be the product of the 0-section in $T^{\ast}B$ with $\mathbb{R}$. Chekanov's theorem \cite{C} states that the number of intersections of $L$ with the 0-wall is greater of equal than $\text{cl}(B)+1$, where $\text{cl}(B)$ denotes the cup-length of $B$. Moreover it is greater or equal than the sum of the Betti numbers of $B$ if we assume that all intersections are transverse. Consider now the contact manifolds $\mathbb{R}^{2n+1}$ and $\mathbb{R}^{2n}\times S^1$. If $M=\mathbb{R}^{2n+1}$ or $M=\mathbb{R}^{2n}\times S^1$ then we have a contact embedding $M\times M\times \mathbb{R}\hookrightarrow J^1M$ sending the diagonal $\Delta$ of $M\times M\times \mathbb{R}$ to the 0-section of $J^1M$ (Bhupal's formula \cite{B}). We denote by $\Gamma_{\phi}$ the Legendrian submanifold of $J^1M$ which is the image of the graph of $\phi$ under this embedding. Then it can be seen that translated points of $\phi$ correspond to Reeb chords between the 0-section and $\Gamma_{\phi}$, i.e. to intersections of $\Gamma_{\phi}$ with the 0-wall. Thus the questions of Problem 1 in the case of $M=\mathbb{R}^{2n+1}$ or $M=\mathbb{R}^{2n}\times S^1$ are answered by applying Chekanov's theorem \cite{C}. We will discuss this in more details in Section \ref{prelim}, taking it as an opportunity to introduce the background material that is needed in the rest of the article. In Sections \ref{per1} and \ref{per2} instead we discuss and prove Problem 2 in the special case of a \textit{positive} contactomorphism of respectively $\mathbb{R}^{2n+1}$ and $\mathbb{R}^{2n}\times S^1$. Recall that a contactomorphism of a closed contact manifold is said to be positive if it is the time-1 map of a positive contact isotopy, i.e. a contact isotopy that is generated by a positive Hamiltonian and thus moves every point in a direction positively transverse to the contact distribution. This notion was introduced by Eliashberg and Polterovich in \cite{EP}. For a compactly supported contactomorphism $\phi$ of $\mathbb{R}^{2n+1}$ or $\mathbb{R}^{2n}\times S^1$, we say that $\phi$ is positive if it is the time-1 map of a Hamiltonian which is positive in the interior of its support. Our main result is the following theorem.

\begin{thm}\label{main}
If $\phi$ is a positive compactly supported contactomorphism of $\mathbb{R}^{2n+1}$ or $\mathbb{R}^{2n}\times S^1$ then $\phi$ has infinitely many geometrically distinct iterated translated points in the interior of its support. 
\end{thm}

As discussed above, Chekanov's theorem and Bhupal's graph construction immediately imply existence of a translated point for every iteration $\phi^k$. As we will see in Section \ref{per1}, a purely formal consequence of this is that either there are infinitely many geometrically distinct iterated translated points of $\phi$, or there is a \textit{periodic point} (i.e. a translated point of some $\phi^k$ which is also a fixed point of $\phi^k$). To get Theorem \ref{main} we need thus to exclude the existence of a periodic point when $\phi$ is a positive contactomorphism. Note that positivity of $\phi$ does not imply that $\phi$ should move every point of  $\mathbb{R}^{2n+1}$ up in the $z$-direction: $\phi$ is the time-1 map of a contact isotopy that moves every point in a direction positively transverse to the contact distribution, but not necessarily in a direction that has a positive scalar product with the Reeb vector field. In fact, it is possible to have positive contact isotopies that move some point always down in the $z$-direction (see \cite{CFP} for concrete examples). To prove Theorem \ref{main} we will use the spectral invariant $c^+$ for contactomorphisms of $\mathbb{R}^{2n+1}$ and $\mathbb{R}^{2n}\times S^1$, that was introduced by Bhupal \cite{B}.\\
\\
Theorem \ref{main} can be seen as a (partial) contact analogue of the result of Viterbo \cite{V} on the existence of infinitely many non-trivial periodic\footnote{In the symplectic case both expressions \textit{periodic point} and \textit{iterated point} are equivalently used to indicate a fixed point of some iteration. Note that if $q$ is a fixed point of some $\phi^k$ then it is also a fixed point of $\phi^{lk}$ for all $l$. Since the equivalent statement is not true for translated points of contactomorphisms, we prefer to use the expression \textit{iterated translated point} instead of \textit{periodic translated point}.} points for a compactly supported Hamiltonian symplectomorphism $\phi$ of $\mathbb{R}^{2n}$. Viterbo's result is proved as follows. Recall that to every compactly supported Hamiltonian symplectomorphism $\phi$ of $\mathbb{R}^{2n}$ we can associate a Lagrangian submanifold $\Gamma_{\phi}$ of $T^{\ast}S^{2n}$ Hamiltonian isotopic to the 0-section. Recall also that to any such Lagrangian submanifold we can associate a generating function quadratic at infinity \cite{Ch1,LS, S, S2}. Fixed points of $\phi$ correspond to intersections of $\Gamma_{\phi}$ with the 0-section and thus to critical points of the generating function $S_{\phi}$ of $\Gamma_{\phi}$. Moreover, critical values of $S_{\phi}$ are given by the symplectic action of the corresponding fixed points of $\phi$. By minimax methods on the generating function we can define spectral invariants $c^+(\phi)$ and $c^-(\phi)$. We have that $c^+(\phi)\geq 0$, $c^-(\phi)\leq 0$, and $c^+(\phi)=c^-(\phi)=0$ if and only if $\phi$ is the identity. For every $k$ we denote by $q_k$ the fixed point of $\phi^k$ corresponding to $c^+(\phi^k)$. We need to show that infinitely many of the $q_k$ for $k\in\mathbb{N}$ are geometrically distinct. Assume first that $\phi$ is non-negative in the sense of Viterbo's partial order, i.e. $c^-(\phi)=0$ and hence $c^+(\phi)>0$. Then $c^+(\phi^k)>0$ for all $k$. The two key facts that allow us to find infinitely many geometrically distinct points among the $q_k$ are the following. 

\begin{enumerate}\label{sympl}
\renewcommand{\labelenumi}{(\roman{enumi})}
\item If $q_{k_1}=q_{k_2}$ then $\frac{1}{k_1}c^+(\phi^{k_1})=\frac{1}{k_2}c^+(\phi^{k_2})$. This is easy to prove using the fact that the symplectic action of a fixed point $q$ of a Hamiltonian symplectomorphism $\phi$ is given by the symplectic area enclosed by $\gamma \sqcup \phi(\gamma)^{-1}$, where $\gamma$ is any path connecting $q$ to some point $q'$ outside the support of $\phi$.
\item The set $\{\,c^+(\phi^k)\;,\; k\in\mathbb{N}_>\,\}$ is bounded above by the capacity of the support of $\phi$.
\end{enumerate}

These two facts together easily imply that infinitely many of the $q_k$ must be geometrically distinct. In the general case, i.e. if $\phi$ is not necessarily non-negative, we can consider a sequence $\phi^{k(i)}$ such that $c^+(\phi^{k(i)})>0$ (we allow $k(i)$ to have any sign). Again (i) and (ii) apply and we get Viterbo's result.\\
\\
As we are going to see in Section \ref{per1}, although we have a contact analogue (introduced by Bhupal \cite{B}) of the spectral invariants $c^+$ and $c^-$, neither (i) nor (ii) hold for contactomorphisms of $\mathbb{R}^{2n+1}$. However in this case we have that if $q$ is a translated point for both $\phi^{k_1}$ and $\phi^{k_2}$ (with $k_1<k_2$) then $\phi^{k_1}(q)$ is a translated point for $\phi^{k_2-k_1}$. This elementary fact, together with monotonicity of $c^+$, allows us to still find infinitely many geometrically distinct iterated translated points for $\phi$. In the case of $\mathbb{R}^{2n}\times S^1$, which is treated in Section \ref{per2}, because of 1-periodicity in the $z$-coordinate the above argument fails for iterated translated points with integer contact action. However in this case an analogue of the rigidity result (ii) holds, and this, together with monotonicity of $c^+$, allows us to conclude.\\
\\
To summarize, the only ingredient in the proof of Theorem \ref{main} in the case of $\mathbb{R}^{2n+1}$ is the existence of a spectral invariant $c^+$ which is carried by translated points and is positive for positive contactomorphisms. Such a spectral invariant was introduced by Bhupal \cite{B}, relying on the work of Chekanov \cite{C}. In the case of $\mathbb{R}^{2n}\times S^1$ we need to use in addition the contact analogue of (ii), i.e. the energy-capacity inequality for spectral invariants that was established in \cite{mio1,mio2}.\\
\\
It is important to notice that the problem of translated points is a special case of the problem of leafwise coisotropic intersections. Let $M$ be a codimension $m$ coisotropic submanifold of a symplectic manifold $(W,\omega)$. The $m$-dimensional distribution $\text{ker}(\omega_{|M})$ is integrable (see for example \cite{MS}) and defines thus a foliation $\mathcal{F}$ on $M$. Given a Hamiltonian symplectomorphism $\phi$ of $W$, a point $q$ of $M$ is called a \textit{leafwise intersection} of $M$ and $\phi(M)$ if $q\in F\cap \phi(F)$ for some leaf $F$ of $\mathcal{F}$. The problem of existence of leafwise intersections was first posed by Moser \cite{M} in 1978. The first results are due to Moser and Banyaga \cite{M,Ba}: they proved existence of leafwise intersections for closed coisotropic submanifolds and Hamiltonian diffeomorphisms $\mathcal{C}^1$-close to the identity. Ekeland and Hofer \cite{EH,H} replaced, for hypersurfaces of restricted contact type in $\mathbb{R}^{2n}$, the $\mathcal{C}^1$-small condition by a much weaker smallness assumption, namely that the Hofer norm of the Hamiltonian symplectomorphism should be smaller than a certain capacity of the region bounded by the hypersurface. This result was then recently extended by Dragnev \cite{D} to coisotropic submanifolds of $\mathbb{R}^{2n}$ with higher codimension, and by Ginzburg \cite{G} to hypersurfaces of restricted contact type in subcritical Stein manifolds. Other results were obtained by Albers, Frauenfelder, McLean and Momin \cite{AF1,AF2,AF3,AF4,AMcL,AM}, Kang \cite{K}, G\"{u}rel \cite{Gu}, Macarini, Merry and Paternain \cite{Me,MMP} and Ziltener \cite{Z}. Translated points of contactomorphisms are a special case of leafwise coisotropic intersections. Indeed, consider the symplectization $\big(M\times \mathbb{R}\,,\,\omega=d(e^{\theta}\alpha)\big)$ of a contact manifold $\big(M\,,\,\xi=\text{ker}(\alpha)\big)$. Then $M$, which we will identify with $M\times\{0\}\subset M\times\mathbb{R}$, is a codimension $1$ coisotropic submanifold of $M\times\mathbb{R}$ and its characteristic foliation coincides with the Reeb foliation. If $\phi$ is a contactomorphism of $M$ which is contact isotopic to the identity then we can lift it to a Hamiltonian symplectomorphism $\widetilde{\phi}$ of $M\times\mathbb{R}$ by setting $\widetilde{\phi}(q,\theta)=\big(\phi(q),\theta-g(q)\big)$ where $g$ is the function defined by $\phi^{\ast}\alpha=e^g\alpha$. Then leafwise intersections for $\widetilde{\phi}$ correspond to translated points of $\phi$. We remark however that none of the existing articles on leafwise coisotropic intersections consider the special case of leafwise intersections corresponding to translated points of contactomorphisms. Notice that, since the question of existence of leafwise intersections is known to be false in general, when studying this problem it is necessary to make some assumption either on the ambient manifold or on the Hamiltonian symplectomorphism. The assumption made so far in the literature is typically a smallness assumption for the Hamiltonian symplectomorphism, which is not suitable to study the problem of leafwise intersections coming from translated points of contactomorphisms. Note however that some of the existing results on leafwise coisotropic intersections  \cite{AF1,AF3,AF4,AMcL,Me,MMP}, \cite[Theorem C]{AF2}, \cite[Theorem 1]{EH} do not need any assumption either on the $\mathcal{C}^1$ or on the Hofer norm of the Hamiltonian symplectomorphism, and therefore it seems possible that they could be applied to obtain existence of translated points of contactomorphisms in some special cases. But, as far as I understand, these results cannot be applied to recover the main theorem of this article.

\subsection*{Acknowledgements} My research was supported by an ANR GETOGA postdoctoral fellowship. I started thinking about the problem of Theorem \ref{main} last year during the Workshop on Symplectic Geometry, Contact Geometry and Interactions in Paris, when Viktor Ginzburg asked me if it was possible to generalize to the contact case Viterbo's result on periodic points of Hamiltonian symplectomorphisms. While the result itself was easy to obtain, it took me a very long time to understand how to put it in the right context. I thank Viktor for listening and commenting on several meaningless versions of Theorem \ref{main} that I told him during the MSRI Program on Symplectic and Contact Geometry and Topology and the Workshop on Symplectic Techniques in Conservative Dynamics in Leiden. And for encouraging me to write it down when I finally got it right, during the Edi-fest in Z\"{u}rich. Moreover I thank him and Peter Albers for their feedback on previous versions of this paper, and the referee for very useful remarks. I also thank the organizers of all the events mentioned above for giving me the opportunity to participate. Finally I thank Vincent Colin, Paolo Ghiggini, Fran\c{c}ois Laudenbach, Marco Mazzucchelli, Josh Sabloff and Lisa Traynor for discussions and support.

\section{Generating functions for contactomorphisms of $\mathbb{R}^{2n+1}$ and $\mathbb{R}^{2n}\times S^1$ and existence of translated points}\label{prelim}

In this section we will give the background information on generating functions for contactomorphisms of $\mathbb{R}^{2n+1}$ and $\mathbb{R}^{2n}\times S^1$ that is needed in the rest of the article. In particular we will recall how Chekanov's theorem and Bhupal's graph construction give existence of generating functions quadratic at infinity, immediately implying the existence of translated points. We refer to \cite{mio1} for more details about the material recalled in this section.\\
\\
Let $\phi$ be a contactomorphism of $\big(\mathbb{R}^{2n+1},\xi=\text{ker}\,(dz-ydx)\big)$ and $g:\mathbb{R}^{2n+1}\rightarrow\mathbb{R}$ the function defined by $\phi^{\ast}(dz-ydx)=e^g(dz-ydx)$. Following Bhupal \cite{B} we define a Legendrian embedding $\Gamma_{\phi}:\mathbb{R}^{2n+1}\longrightarrow J^1\mathbb{R}^{2n+1}$ to be the composition $\tau \circ \text{gr}_{\phi}$, where $\text{gr}_{\phi}:\mathbb{R}^{2n+1}\longrightarrow \mathbb{R}^{2(2n+1)+1}$ is the Legendrian embedding $q\mapsto (q,\phi(q),g(q))$ and $\tau:\mathbb{R}^{2(2n+1)+1}\longrightarrow J^1\mathbb{R}^{2n+1}$ the contact embedding $(x,y,z,X,Y,Z,\theta)\mapsto \big(x,Y,z, Y-e^{\theta}y, x-X, e^{\theta}-1, xY-XY+Z-z\big)$. Here we consider the product contact structure $e^{\theta}(dz-ydx)-(dZ-YdX)$ on $\mathbb{R}^{2(2n+1)+1}$ and the standard contact structure $d\theta-pdq$ on $J^1\mathbb{R}^{2n+1}$. More explicitly, for $\phi=(\phi_1,\phi_2,\phi_3)$ we have that $\Gamma_{\phi}:\mathbb{R}^{2n+1}\longrightarrow J^1\mathbb{R}^{2n+1}$ is given by
$$
\Gamma_{\phi}(x,y,z)=(x,\phi_2,z,\phi_2-e^gy,x-\phi_1,e^g-1,x\phi_2-\phi_1\phi_2+\phi_3-z).
$$
If $\phi$ is contact isotopic to the identity then $\Gamma_{\phi}$ is contact isotopic to the 0-section. Moreover, if $\phi$ is compactly supported then we can see $\Gamma_{\phi}$ as a Legendrian submanifold of $J^1S^{2n+1}$.  Similarly, if $\phi$ is a compactly supported contactomorphism of $\mathbb{R}^{2n}\times S^1$ (seen as a contactomorphism of $\mathbb{R}^{2n+1}$ which is 1-periodic in the $z$-coordinate) then we can see $\Gamma_{\phi}$ as a Legendrian submanifold of $J^1\big(S^{2n}\times S^1\big)$. \\
\\
Notice that, since the Reeb vector field in $\big(\mathbb{R}^{2n+1},\xi=\text{ker}\,(dz-ydx)\big)$ is given by $\frac{\partial}{\partial z}$, translated points of $\phi$ are points $q=(x,y,z)$ with $\phi_1(q)=x$, $\phi_2(q)=y$ and $g(q)=0$. They correspond to intersections of $\Gamma_{\phi}$ with the 0-wall and so to Reeb chords between $\Gamma_{\phi}$ and the 0-section. We are thus in the setting of Chekanov's theorem on critical points of quasi-functions \cite{C}, that we will now recall.\\
\\
Let $B$ be a smooth closed manifold, and consider the 1-jet bundle $J^1B=T^{\ast}B\times \mathbb{R}$ equipped with the contact form $d\theta-pdq$. Note that the 1-jet $j^1f=\{\;\big(x,df(x),f(x)\big)\,|\,x\in B\;\}$ of a smooth function $f:B\rightarrow \mathbb{R}$ is a Legendrian submanifold of $J^1B$ and that (transverse) intersections of $j^1f$ with the 0-wall, which is the product in $J^1B=T^{\ast}B\times \mathbb{R}$ of the 0-section with $\mathbb{R}$, correspond to (non-degenerate) critical points of $f$. Hence we can give a lower bound on the number of intersections of $j^1f$ with the 0-wall in terms of the topology of $B$, by applying the Lusternik-Shnirelman and (in the non-degenerate case) Morse inequalities. Chekanov proved that these estimates remain valid also for arbitrary Legendrian submanifolds of $J^1B$ that are contact isotopic to the 0-section (he calls such Legendrian submanifolds \textit{quasi-functions} on $B$). With this he generalized the analogous result by Laudenbach and Sikorav \cite{LS} on the existence of intersections with the 0-section of Lagrangian submanifolds of $T^{\ast}B$ Hamiltonian isotopic to the 0-section. Chekanov's theorem follows from the existence of a generating function quadratic at infinity for every Legendrian submanifold of $J^1B$ contact isotopic to the 0-section, and from Conley-Zehnder theorem on critical points of functions quadratic at infinity \cite{CZ}.\\
\\
Recall that a function $S:E\rightarrow\mathbb{R}$, where $E$ is the total space of a fiber bundle $\pi:E \longrightarrow B$, is called a \textit{generating function} for a Legendrian submanifold $L$ of $J^1B$ if $dS: E\longrightarrow T^{\ast}E$ is transverse to $N_E:=\{\,(e,\eta)\in T^{\ast}E \;|\; \eta = 0 \;\text{on} \;\text{ker}\,d\pi\,(e)\,\}$ and if $L$ is the image of the Legendrian immersion $j_S:\Sigma_S\rightarrow J^1B$, $e\mapsto\big(\pi(e), v^{\ast}(e),S(e)\big)$. Here $\Sigma_S\subset E$ is the set of fiber critical points of $S$ and for $e\in\Sigma_S$ the element $v^{\ast}(e)$ of $T^{\phantom{\pi}\ast}_{\pi(e)}B$ is defined by $v^{\ast}(e)\,(X):=dS\,(\widehat{X})$ for $X \in T_{\pi(e)}B$, where $\widehat{X}$ is any vector in $T_eE$ with $\pi_{\ast}(\widehat{X})=X$. Note that critical points of $S$ correspond under $j_S$ to intersections of $L$ with the 0-wall. Moreover we have that non-degenerate critical points correspond to transverse intersections. Note also that if $e$ is a critical point of $S$ then its critical value is given by the $\mathbb{R}$-coordinate of the point $j_S(e)$. A generating function $S:E \longrightarrow \mathbb{R}$ is said to be \textit{quadratic at infinity} if $\pi:E \longrightarrow B $ is a vector bundle and if there exists a non-degenerate quadratic form $Q_{\infty}: E \longrightarrow \mathbb{R}$ such that $dS-\partial_vQ_{\infty}: E \longrightarrow E^{\ast}$ is bounded, where $\partial_v$ denotes the fiber derivative. As already mentioned, existence of generating functions quadratic at infinity for Legendrian submanifolds of $J^1B$ contact isotopic to the 0-section was proved by Chekanov \cite{C}. See also Chaperon \cite{Ch} for a different proof of the same result.\\
\\
Given a contactomorphism $\phi$ of $\mathbb{R}^{2n+1}$ or $\mathbb{R}^{2n}\times S^1$, compactly supported and isotopic to the identity, by applying Chekanov's theorem to the corresponding Legendrian submanifold $\Gamma_{\phi}$ of $J^1S^{2n+1}$ or $J^1(S^{2n}\times S^1)$ we get thus a lower bound on the number of translated points. In the case of $\mathbb{R}^{2n+1}$ we get at least one translated point in the interior of the support (the generating function has at least two critical points, but one of them might correspond to the point at infinity of $\mathbb{R}^{2n+1}$). Similarly in the case of $\mathbb{R}^{2n}\times S^1$ we get at least two translated points in the interior of the support, and at least three if all of them are non-degenerate (i.e. if they correspond to transverse intersections of $\Gamma_{\phi}$ with the 0-wall).\\
\\ 
In the next sections we will prove that every positive contactomorphism $\phi$ of $\mathbb{R}^{2n+1}$ or $\mathbb{R}^{2n}\times S^1$ which is compactly supported and isotopic to the identity has infinitely many iterated translated points in the interior of the support. By the discussion above we know that every iteration $\phi^k$ has at least one translated point in the interior of the support. However these iterated translated points are not necessarily geometrically distinct. To find infinitely many geometrically distinct iterated translated points out of them we will use the spectral invariants $c^{\pm}(\phi)$ that were introduced in \cite{B}. Recall that the definition of $c^{\pm}(\phi)$ is based on the following minimax method for generating functions \cite{V}. Let $L$ be a Legendrian submanifold of $J^1B$ with generating function $S:E\rightarrow \mathbb{R}$. Denote by $E^{a}$, for $a\in \mathbb{R}\cup\infty$, the sublevel set of $S$ at $a$, and by $E^{-\infty}$ the set $E^{-a}$ for $a>0$ big. We consider the inclusion $i_a: (E^a,E^{-\infty})\hookrightarrow (E,E^{-\infty})$, and the induced map on cohomology
$$ i_a^{\phantom{a}\ast}: H^{\ast}(B)\equiv H^{\ast}(E,E^{-\infty})\longrightarrow H^{\ast}(E^a,E^{-\infty})$$
where $H^{\ast}(B)$ is identified with $H^{\ast}(E,E^{-\infty})$ via the Thom isomorphism. For any $u\neq 0$ in $H^{\ast}(B)$ we define
$$ c(u,L)=c(u,S)=\text{inf}\,\{\,a\in \mathbb{R} \;|\;i_a^{\phantom{a}\ast}(u)\neq 0\,\}.$$
Note that $c(u,L)$ is a critical value of $S$. The good definition of $c(u,L)$ follows from the uniqueness theorem for generating functions quadratic at infinity \cite{V,Th,Th2}. For a contactomorphism $\phi$ of $\mathbb{R}^{2n+1}$ or $\mathbb{R}^{2n}\times S^1$ we define $c^+(\phi):=c(\mu,\Gamma_{\phi})$ and $c^-(\phi):=c(1,\Gamma_{\phi})$ where $\mu$ and $1$ are respectively the orientation and the unit cohomology class either of $S^{2n+1}$ or $S^{2n}\times S^1$. It can be proved that $c^+(\phi)\geq 0$, $c^-(\phi)\leq 0$, and $c^+(\phi)=c^-(\phi)=0$ if and only if $\phi$ is the identity.\\
\\
If $q=(x,y,z)$ is a translated point of $\phi$ then we call the real number $\mathcal{A}_{\phi}(q):=\phi_3(q)-z$ its \textit{contact action} (in the case of $\mathbb{R}^{2n}\times S^1$ we regard $\phi$ as a 1-periodic contactomorphism of $\mathbb{R}^{2n+1}$ so that also in this case $\phi_3(q)-z$ is a well-defined real number). Note that $\mathcal{A}_{\phi}(q)$ is equal to the critical value of the critical point of the generating function of $\phi$ corresponding to $q$. It follows from this discussion that for every contactomorphism $\phi$ of $\mathbb{R}^{2n+1}$ or $\mathbb{R}^{2n}\times S^1$ there are two translated points $q^+$ and $q^-$ such that $\mathcal{A}_{\phi}(q^+)=c^+(\phi)$ and $\mathcal{A}_{\phi}(q^-)=c^-(\phi)$. Unless $\phi$ is the identity, at least one of them has non-trivial contact action, and so in particular belongs to the interior of the support of $\phi$.

\section{Iterated translated points in $\mathbb{R}^{2n+1}$}\label{per1}

Let $\phi$ be a contactomorphism of $\mathbb{R}^{2n+1}$, compactly supported and isotopic to the identity (but different than the identity). By the discussion in Section \ref{prelim} we know that $\phi$ has at least one translated point in the interior of its support. Since the same is true for all the iterations $\phi^k$, we have infinitely many iterated translated points of $\phi$ in the interior of the support. The question is whether infinitely many of them are geometrically distinct. In the following theorem we will study this problem in the case of a positive contactomorphism $\phi$. Recall that $\phi$ is called positive if it is the time-1 flow of a contact Hamiltonian which is positive in the interior of its support. Thus, $\phi$ is the time-1 map of a contact isotopy that in the interior of its support moves every point in a direction positively transverse to the contact distribution. As it was proved by Bhupal \cite{B} this relation induces a well-defined bi-invariant partial order on the group of compactly supported contactomorphisms of $\mathbb{R}^{2n+1}$ isotopic to the identity.

\begin{thm}\label{R}
A positive compactly supported contactomorphism $\phi$ of $\mathbb{R}^{2n+1}$ has infinitely many non-trivial geometrically distinct iterated translated points.
\end{thm}

Here an iterated translated point of $\phi$ is considered non-trivial if it has a non-zero contact action, and so in particular belongs to the interior of the support of $\phi$.\\
\\
To prove Theorem \ref{R} we need the following lemma.

\begin{lemma}\label{transl_obvious}
Let $k_1$ and $k_2$ be positive integers with $k_1<k_2$, and assume that $q$ is a translated point for both $\phi^{k_1}$ and $\phi^{k_2}$. Then $\phi^{k_1}(q)$ is a translated point for $\phi^{k_2-k_1}$, hence in particular an iterated translated point for $\phi$.
\end{lemma}

\begin{proof}
Since $q$ is a translated point for $\phi^{k_1}$ and $\phi^{k_2}$ we have in particular that $\phi^{k_1}(q)$ and $\phi^{k_2}(q)=\phi^{k_2-k_1}\big(\phi^{k_1}(q)\big)$ are in the same Reeb orbit. For any positive integer $k$ we denote by $g_k$ the function satisfying $(\phi^k)^{\ast}\alpha=e^{g_k}\alpha$. Note that $g_k=g\circ \phi^{k-1}+g\circ \phi^{k-2}+\cdots + g\circ \phi + g$, and so in particular $g_{k_2}(q)=g_{k_2-k_1}\big(\phi^{k_1}(q)\big)+g_{k_1}(q)$. Since $g_{k_1}(q)=0=g_{k_2}(q)$ we get that $g_{k_2-k_1}\big(\phi^{k_1}(q)\big)=0$, concluding the proof that $\phi^{k_1}(q)$ is a translated point for $\phi^{k_2-k_1}$.
\end{proof}

The key ingredient in the proof of Theorem \ref{R} is strict monotonicity of $c^+$: if $\phi$ and $\psi$ are positive contactomorphisms with $\phi<\psi$ (i.e. $\psi\phi^{-1}$ is positive) then $c^+(\phi)<c^+(\psi)$. Note that the proof of monotonicity of $c^{\pm}$ in \cite{mio1} also gives strict monotonicity. In particular, if $\phi$ is positive then we have that $0<c^+(\phi)<c^+(\phi^2)<c^+(\phi^3)<\cdots$.

\begin{proof}[Proof of Theorem \ref{R}]
For every positive integer $k$ we denote by $q_k$ the translated point of $\phi^k$ whose contact action is equal to $c^+(\phi^k)$ (see the discussion at the end of the previous section). Then either infinitely many of the $\{\,q_k\,,\, k\in\mathbb{N}_>\,\}$ are geometrically distinct or there exist a point $q$ and an infinite sequence $k_1<k_2<k_3<\cdots$ such that $q=q_{k_i}$ for all $i$ in $\mathbb{N}_>$. But in the second case we can consider the points $\{\,\phi^{k_i}(q)\,,\, i\in\mathbb{N}_>\,\}$, that by Lemma \ref{transl_obvious} are iterated translated points of $\phi$. These points are all geometrically distinct. Indeed, let $q=(x,y,z)$ and $\phi^{k_i}(q)=(x,y,z_{k_i})$. Then $c^+(\phi_{k_i})=z_{k_i}-z$ and thus by monotonicity of $c^+$ and positivity of $\phi$ we get that $z_{k_1}<z_{k_2}<z_{k_3}<\cdots$ so that $\phi_{k_1}(q)$, $\phi_{k_2}(q)$, $\phi_{k_3}(q)$, $\cdots$ are all distinct.
\end{proof}

We would like to stress that the fact that, in the previous proof, the points $\{\,\phi^{k_i}(q)\,,\, i\in\mathbb{N}_>\,\}$ are distinct is not obvious a priori and really needs some deep tool to be proved, in our case monotonicity of the spectral invariant $c^+$. Indeed, as discussed in the introduction, the positivity assumption for a contactomorphism $\phi$ of $\mathbb{R}^{2n+1}$ does not imply that this contactomorphism should move every point up in the $z$-direction. \\
\\
If we remove the positivity assumption in Theorem \ref{main}, we can just say that a compactly supported contactomorphism $\phi$ of $\mathbb{R}^{2n+1}$ either has a \textit{periodic point} (i.e. a fixed point $q$ of some $\phi^k$, with $g_k(q)=0$) or infinitely many non-trivial geometrically distinct iterated translated points. Indeed, suppose as in the proof of Theorem \ref{R} that there exist a point $q$ and an infinite sequence $k_1<k_2<k_3<\cdots$ such that $q=q_{k_i}$ for all $i$ in $\mathbb{N}_>$. Then either the set $\{\,\phi^{k_i}(q)\,,\, i\in\mathbb{N}_>\,\}$ has infinitely many distinct points or there exist $k_1<k_2$ such that $\phi^{k_1}(q)=\phi^{k_2}(q)=:q'$. But then, by Lemma \ref{transl_obvious}, $q'$ is a fixed point of $\phi^{k_2-k_1}$ with $g_{k_2-k_1}=0$. Notice that this argument is only formal, and would apply to any contact manifold once existence of a non-trivial translated point for each iterate $\phi^k$ of a contactomorphism $\phi$ is established.\\
\\
In the study of iterated translated points for contactomorphisms of $\mathbb{R}^{2n+1}$, the main difference with respect to the symplectic case is that we do not have any rigidity of the contact action (i.e. something like properties (i) and (ii) on page \pageref{sympl}). The only tool we have in this case to distinguish the iterated translated points is monotonicity of $c^+$. As we are going to see in the next section, in $\mathbb{R}^{2n}\times S^1$ the automatic distinction of iterated translated points given by monotonicity of $c^+$ fails for integer contact actions because of 1-periodicity in the $z$-coordinate. On the other hand in $\mathbb{R}^{2n}\times S^1$ we do have some rigidity results for the contact action, that will allow us to still find infinitely many iterated translated points.

\section{Iterated translated points in $\mathbb{R}^{2n}\times S^1$}\label{per2}

We will show in this section the existence of infinitely many geometrically distinct iterated translated points for positive compactly supported contactomorphisms of $\mathbb{R}^{2n}\times S^1$. The argument given in the case of $\mathbb{R}^{2n+1}$ now fails since, because of 1-periodicity in the $z$-coordinate, monotonicity of $c^+$ is not enough any more to distinguish, in the case $q_{k_1}=q_{k_2}$, the points $\phi^{k_1}(q_{k_1})$ and $\phi^{k_2}(q_{k_2})$. Indeed, if $c^+(\phi^{k_2})-c^+(\phi^{k_1})$ is integer then $\phi^{k_1}(q_{k_1})=\phi^{k_2}(q_{k_2})$. However we still get existence of infinitely many iterated translated points thanks to some rigidity results for the contact action. In particular we have the following new phenomena with respect to the case of $\mathbb{R}^{2n+1}$.\\
\\
For contactomorphisms of $\mathbb{R}^{2n}\times S^1$ we have that the integer parts of $c^+$ and $c^-$ are invariant by conjugation, i.e. $\lceil c^{\pm}(\phi)\rceil=\lceil c^{\pm}(\psi\phi\psi^{-1})\rceil$ and $\lfloor c^{\pm}(\phi)\rfloor=\lfloor c^{\pm}(\psi\phi\psi^{-1})\rfloor$. Moreover the following triangle inequalities hold: $$\lceil c^+(\phi\psi)\rceil\leq\lceil c^+(\phi)\rceil+\lceil c^+(\psi)\rceil$$ 
$$\lfloor c^-(\phi\psi)\rfloor\geq\lfloor c^-(\phi)\rfloor+\lfloor c^-(\psi)\rfloor.$$
These properties, that follow from 1-periodicity in the $z$-coordinate (see \cite{mio1}), allow us to define a bi-invariant metric on the contactomorphism group of $\mathbb{R}^{2n}\times S^1$ by
$$d\,(\phi,\psi):=\lceil c^+(\phi\psi^{-1})\rceil-\lfloor c^-(\phi\psi^{-1})\rfloor$$
and a contact capacity for a domain $\mathcal{V}$ of $\mathbb{R}^{2n}\times S^1$ by 
$$c(\mathcal{V}):=\text{sup}\,\{\,\lceil c^+(\phi)\rceil \:|\: \phi\in\text{Cont}\,(\mathcal{V})\,\}$$ 
where $\text{Cont}\,(\mathcal{V})$ denotes the set of time-1 maps of contact Hamiltonians supported in $\mathcal{V}$. The contact capacity $c(\mathcal{V})$ is well-defined because it can be proved that if $\psi$ is a contactomorphism of $\mathbb{R}^{2n}\times S^1$ displacing $\mathcal{V}$, i.e. such that $\psi(\mathcal{V})\cap\mathcal{V}= \emptyset$, then $\lceil c^+(\phi)\rceil \leq E(\psi)$ for every $\phi$ in $\text{Cont}\,(\mathcal{V})$. Here $E(\psi)$ is the energy of $\psi$ with respect to the metric $d$, i.e. $E(\phi):=\lceil c^+(\phi)\rceil-\lfloor c^-(\phi)\rfloor$. \\
\\
We refer to \cite{mio1,mio2} for a proof and further discussion of these results. To study iterated translated points of a contactomorphism $\phi$ of $\mathbb{R}^{2n}\times S^1$ we are only going to use the fact that the sequence $\{\,c^+(\phi^k)\,,\, k\in\mathbb{N}_>\,\}$ is bounded above, by the capacity of the support of $\phi$. This, together with monotonicity of $c^+$, allows us to prove the following theorem.

\begin{thm}\label{existence_S}
If $\phi$ is a positive contactomorphism of $\mathbb{R}^{2n}\times S^1$ which is compactly supported and isotopic to the identity then $\phi$ has infinitely many non-trivial geometrically distinct iterated translated points.
\end{thm}

\begin{proof}
For every positive integer $k$ we consider the translated point $q_k$ of $\phi^k$ whose contact action is equal to $c^+(\phi^k)$. As in the proof of Theorem \ref{R}, either infinitely many of the $\{\,q_k\,,\, k\in\mathbb{N}_>\,\}$ are geometrically distinct or there exist a point $q$ and an infinite sequence $k_1<k_2<k_3<\cdots$ such that $q=q_{k_i}$ for all $i$ in $\mathbb{N}_>$. In this second case we can consider the points $\{\,\phi^{k_i}(q)\,,\, i\in\mathbb{N}_>\,\}$. Note that by Lemma \ref{transl_obvious} all of them are iterated translated points of $\phi$. Let $q=(x,y,z)$ and $\phi^{k_i}(q)=(x,y,z_{k_i})$. Because of monotonicity of $c^+$, either infinitely many of the $\{\,\phi^{k_i}(q)\,,\, i\in\mathbb{N}_>\,\}$ are geometrically distinct or there is an infinite subset $I$ of $\mathbb{N}_>$ such that $z_{k_i}-z\in\mathbb{Z}$ for all $i$ in $I$. But this second possibility gives a contradiction because the sequence $\{\,c^+(\phi^{k_i})=z_{k_i}-z\,,\, i\in I\,\}$ is increasing, because of monotonicity of $c^+$ and positivity of $\phi$, and bounded above by the capacity of the support of $\phi$ (see the discussion above) and thus cannot be made of integers.
\end{proof}


\begin{thebibliography}{99}

\bibitem[Abb99a]{Ab1} C.~Abbas, A note on V.I. Arnold's Chord Conjecture, \textit{Int. Math. Res. Notices} \textbf{4} (1999), 217--222.

\bibitem[Abb99b]{Ab2} C.~Abbas, Finite energy surfaces and the Chord Problem,  \textit{Duke Math. J.} \textbf{96} (1999), 241--316.

\bibitem[AF08]{AF1} P.~Albers and U.~Frauenfelder, Infinitely many leaf-wise intersections on cotangent bundles, \textit{arXiv}: 0812.4426.

\bibitem[AF10a]{AF2} P.~Albers and U.~Frauenfelder, Leaf-wise intersections and Rabinowitz Floer homology, \textit{J. Topol. Anal.} \textbf{2} (2010), 77--98.

\bibitem[AF10b]{AF3} P.~Albers and U.~Frauenfelder, Spectral invariants in Rabinowitz-Floer homology and global Hamiltonian perturbations, \textit{J. Mod. Dyn.} \textbf{4} (2010), 329--357.

\bibitem[AF10c]{AF4} P.~Albers and U.~Frauenfelder, On a Theorem by Ekeland-Hofer, \textit{arXiv}:1001.3386. 

\bibitem[AMcL09]{AMcL} P.~Albers and M.~McLean, Non-displaceable contact embeddings and infinitely many leaf-wise intersections, \textit{arXiv}:0904.3564.

\bibitem[AM10]{AM} P.~Albers and A.~Momin, Cup-length estimates for leaf-wise intersections, \textit{Math. Proc. Cambridge Philos. Soc.} \textbf{149} (2010), 539--551

\bibitem[Arn86]{A} V.~I.~Arnold, First steps in symplectic topology, \textit{Russ. Math. Surv.} \textbf{41} (1986), 1--21.

\bibitem[Ban80]{Ba} A.~Banyaga, On fixed points of symplectic maps, \textit{Invent. Math.} \textbf{56}(1980), 215--229.

\bibitem[Bh01]{B}
M.~Bhupal, A partial order on the group of contactomorphisms of $\mathbb{R}^{2n+1}$ via generating functions, \textit{Turkish J. Math.} \textbf{25} (2001), 125--135.

\bibitem[Chap84]{Ch1}
M.~Chaperon, Une id\'{e}e du type \textquotedblleft g\'{e}od\'{e}siques bris\'{e}es\textquotedblright  pour les syst\'{e}mes hamiltoniens, \textit{C. R. Acad. Sci. Paris}, S\'{e}r. I Math. \textbf{298} (1984), 293--296. 

\bibitem[Chap95]{Ch} 
M.~Chaperon, On generating families, in \textit{The Floer Memorial Volume} (H. Hofer et al., eds.),
(Progr. Math., vol. 133) Birkhauser, Basel 1995, pp. 283--296.

\bibitem[Chek96]{C}
Y.~Chekanov, Critical points of quasi-functions and generating families of Legendrian manifolds, \textit{Funct. Anal. Appl.} \textbf{30} (1996), 118--128.

\bibitem[Cie02]{Ci} K.~Cieliebak, Handle attaching in symplectic homology and the chord conjecture, \textit{J. Eur. Math. Soc.} 
\textbf{4} (2002), 115--142. 

\bibitem[CFP10]{CFP} V.~Colin, E.~Ferrand and P.~Pushkar, Positive isotopies of Legendrian submanifolds and applications, \textit{arXiv:1004.5263}.

\bibitem[CZ83]{CZ} C.~Conley and E.~Zehnder, The Birkhoff-Lewis fixed point theorem and a conjecture of V.~I.~Arnol'd, \textit{Invent. Math.} \textbf{73} (1983), 33--49. 

\bibitem[Dra08]{D}
D.~L.~Dragnev, Symplectic rigidity, symplectic fixed points, and global perturbations of Hamiltonian systems, \textit{Comm. Pure Appl. Math.} \textbf{61} (2008), 346--370. 

\bibitem[EH89]{EH} I.~Ekeland and H.~Hofer, Two symplectic fixed-point theorems with applications to Hamiltonian dynamics, \textit{J. Math. Pures Appl.} (9) \textbf{68} (1989), no 4, 467--489 (1990).

\bibitem[EP00]{EP}
Y.~Eliashberg and L.~Polterovich, Partially ordered groups and geometry of contact transformations, \textit{Geom. Funct. Anal.} \textbf{10} (2000), 1448--1476.

\bibitem[Gin07]{G} V.~Ginzburg, Coisotropic intersections, \textit{Duke Math. J.} \textbf{140} (2007), 111--163.


\bibitem[Giv90]{Giv} A.~Givental, Nonlinear generalization of the Maslov index, in \textit{Theory
of singularities and its applications}, pp. 71--103, Adv. Soviet Math., 1, Amer. Math. Soc., Providence, RI, 1990.

\bibitem[Gur10]{Gu} B.~G\"{u}rel, Leafwise coisotropic intersections, \textit{Int. Math. Res. Not.} \textbf{2010} 914--931. 

\bibitem[Hof90]{H} H.~Hofer, On the topological properties of symplectic maps,  \textit{Proc. Roy. Soc. Edinburgh}  \textbf{115}  (1990), 25--38.

\bibitem[HT10]{HT} M.~Hutchings and C.~H.~Taubes, Proof of the Arnold chord conjecture in three dimensions I, \textit{arXiv}: 1004.4319. to appear in \textit{Mathematical Research Letters}.

\bibitem[Kan09]{K} J.~Kang, Existence of leafwise intersection points in the unrestricted case, \textit{arXiv}: 09102369.

\bibitem[LS85]{LS} F.~Laudenbach and J.C.~Sikorav, Persistance d'intersection avec la section nulle au cours d'une isotopie hamiltonienne dans un fibre cotangent, \textit{Invent. Math.} \textbf{82} (1985), 349--357.

\bibitem[MMP11]{MMP} L.~Macarini, W.~Merry and G.~P.~Paternain, On the growth rate of leaf-wise intersections, \textit{arXiv}: 1101.4812.

\bibitem[MS]{MS}
D.~McDuff, D.~Salamon, \textit{Introduction to Symplectic Topology}, Oxford University Press, 1998.

\bibitem[Mer10]{Me} W.~Merry, On the Rabinowitz Floer homology of twisted cotangent bundles, \textit{arXiv}: 1002.0162.

\bibitem[Moh01]{Mo} K.~Mohnke, Holomorphic disks and the chord conjecture, \textit{Ann. of Math.}(2) \textbf{154} (2001), 219--222. 

\bibitem[Mos78]{M}
J.~Moser, A fixed point theorem in symplectic geometry, \textit{Acta Math.} \textbf{141}(1978), 17--34.

\bibitem[S11]{mio1}
S.~Sandon, Contact Homology, Capacity and Non-Squeezing in $\mathbb{R}^{2n}\times S^{1}$ via Generating Functions, \textit{Ann. Inst. Fourier (Grenoble)} \textbf{61} (2011), 145 -- 185.

\bibitem[S10]{mio2}
S.~Sandon, An integer valued bi-invariant metric on the group of contactomorphisms of $\mathbb{R}^{2n}\times S^{1}$, \textit{J. Topol. Anal.} \textbf{2} (2010), 327--339.

\bibitem[Sik86]{S}
J.C.~Sikorav, Sur les immersions lagrangiennes dans un fibr\'{e} cotangent admettant une phase g\'{e}n\'{e}ratrice globale, \textit{C.R. Acad. Sci. Paris}, S\'{e}r. I Math. \textbf{302} (1986), 119--122.

\bibitem[Sik87]{S2}
J.C.~Sikorav, Problemes d'intersections et de points fixes en g\'{e}om\'{e}trie hamiltonienne, \textit{Comment. Math. Helv.} \textbf{62} (1987), 62--73.

\bibitem[Th95]{Th}
D.~Th\'{e}ret, Utilisation des fonctions g\'{e}n\'{e}ratrices en g\'{e}om\'{e}trie symplectique globale, Ph.D. Thesis, Universit\'{e} Denis Diderot (Paris 7), 1995.

\bibitem[Th99]{Th2}
D.~Th\'{e}ret, A complete proof of Viterbo's uniqueness theorem on generating functions, \textit{Topology Appl.} \textbf{96} (1999), 249-266.

\bibitem[Vit92]{V}
C.~Viterbo, Symplectic topology as the geometry of generating functions, \textit{Math. Ann.} \textbf{292} (1992), 685--710.

\bibitem[Zil10]{Z} F.~Ziltener, Coisotropic submanifolds, leaf-wise fixed points, and presymplectic embeddings, \textit{J. Symplectic Geom.} \textbf{8}(2010), 95--118.

\end{thebibliography}
\end{document}